\documentclass[10pt]{amsart}
\usepackage{amsmath,amssymb,amsthm,latexsym}

\newtheorem{theorem}{Theorem}

\newtheorem{cor}{Corollary}[theorem]
\newtheorem{conj}[theorem]{Conjecture}

\DeclareMathOperator{\E}{{\mathbb{E}}}

\DeclareMathOperator{\CT}{CT}  

\newcommand{\R}{{\mathbb{R}}}

\renewcommand{\Re}{{\mathfrak{Re}}}

\renewcommand{\O}{{\mathcal{O}}}

\renewcommand{\a}{\alpha}
\renewcommand{\b}{\beta}
\newcommand{\g}{\gamma}

\renewcommand{\t}{\theta}

\renewcommand{\idotsint}{\int\!\cdots\int}
\renewcommand{\i}{{\mathrm{i}}}
\renewcommand{\d}{{\mathrm{d}}}
\newcommand{\Hyper}[3]{{}_2F_1\left(\genfrac{}{}{0pt}{}{#1}{#2};#3\right)}

\numberwithin{equation}{section}

\begin{document}
\title[Discrete moments of the zeta function]{Random matrix theory and discrete moments of the Riemann zeta function}
\author{C.P. Hughes}
\address{Raymond and Beverly Sackler School of Mathematical Sciences,
Tel Aviv University, Tel Aviv 69978, Israel. Current address:
American Institute of Mathematics, 360 Portage Ave., Palo Alto, CA
94306-2244 ({\tt hughes@aimath.org})}

\date{18 November 2002}

\begin{abstract}
We calculate the discrete moments of the
characteristic polynomial of a random unitary matrix, evaluated a small distance away from an
eigenangle. Such results allow us to make conjectures about similar
moments for the Riemann zeta function, and provide a uniform approach
to understanding moments of the zeta function and its derivative.
\end{abstract}

\maketitle

\section{Introduction}\label{sect:1}

Let $\zeta(s)$ be the Riemann zeta function, and denote its
non-trivial zeros by $\tfrac{1}{2}+\i\gamma_n$, with $0<\g_1\leq
\g_2\leq \dots$. (For simplicity, we will assume the Riemann Hypothesis, which says that $\g_n\in\R$). It is known (see, for example, Titchmarsh's book~\cite{Tit1} for details) that if $N(T)$ is the number of zeros with $0<\g_n\leq T$ then
\begin{equation*}
N(T) = \frac{T}{2\pi}\log\frac{T}{2\pi e} + \O(\log T) .
\end{equation*}

Moments of the zeta function,
\begin{equation*}
I_k(T) := \frac{1}{T}\int_0^T \left|\zeta(\tfrac{1}{2}+\i t)\right|^{2k}\;\d t ,
\end{equation*}
have long been of interest to number theorists, with it being widely believed that
\begin{equation*}
\frac{1}{T}\int_0^T \left|\zeta(\tfrac{1}{2}+\i t)\right|^{2k}\;\d t \sim f_k a(k) (\log T)^{k^2}
\end{equation*}
with
\begin{equation}\label{eq:a(k)}
a(k) = \prod_{\substack{p\ \mathrm{prime}}} \left(1-\frac{1}{p}\right)^{k^2} \sum_{m=0}^\infty \left(\frac{\Gamma(m+k)}{m!\,\Gamma(k)}\right)^2p^{-m} ,
\end{equation}
and $f_k$ being an integer (when $k$ is integer) divided by $(k^2)!$, whose value was unknown apart from in a few cases. The known values of $f_k$ are $f_0=1$ (trivial), $f_1=1$ (Hardy and Littlewood~\cite{HardyLittle}), and $f_2=\frac{1}{12}$ (Ingham~\cite{Ingham}). The value $f_3=\frac{42}{9!}$ has been conjectured by Conrey and Ghosh~\cite{ConGhosh2} and $f_4=\frac{24024}{16!}$ is a conjecture of Conrey and Gonek~\cite{ConGonek1}.
We should mention that $a(k)$ given in \eqref{eq:a(k)} can be calculated for certain $k$: $a(0) = a(1) = 1$ and $a(-1)=a(2) = \frac{6}{\pi^2}$.

In~\cite{KS1} Keating and Snaith argued that one can create a
probabilistic model for the zeta function around height $T$ using the
characteristic polynomial of an $N\times N$ unitary matrix chosen
according to Haar measure, when
\begin{equation*}
N=\log\frac{T}{2\pi} .
\end{equation*}

Setting
\begin{align*}
Z_U(\t) &:= \det(I-e^{-\i\t}U)\\
&=\prod_{n=1}^N \left(1-e^{\i(\t_n-\t)}\right) ,
\end{align*}
and defining $M_N(2k):=\E_N\left\{|Z_U(0)|^{2k}\right\}$ where $\E_N$ denotes expectation with respect to Haar measure, they found that
\begin{align}\label{eq:defn M_N(2k)}
M_N(2k) &= \prod_{j=1}^N \frac{\Gamma(j)\Gamma(j+2k)}{\Gamma(j+k)^2}\\
&\sim \frac{G^2(k+1)}{G(2k+1)} N^{k^2} \nonumber
\end{align}
as $N\to\infty$ for fixed $k$ subject to $\Re(k)>-1/2$, where $G(\cdot)$ is the Barnes $G$--function.

By comparing with the known (and previously conjectured) values of
$f_k$, they were led to conjecture that $f_k=\frac{G^2(k+1)}{G(2k+1)}$.
\begin{conj}\label{conj:KS}
{\bf (Keating and Snaith).}
For fixed $k>-1/2$, as $T\to\infty$,
\begin{equation*}
\frac{1}{T}\int_0^T \left|\zeta(\tfrac{1}{2}+\i t)\right|^{2k} \;\d t \sim \frac{G^2(k+1)}{G(2k+1)} a(k) \left(\log \frac{T}{2\pi}\right)^{k^2}
\end{equation*}
where $G(\cdot)$ is the Barnes $G$--function and $a(k)$ is given in \eqref{eq:a(k)}
\end{conj}

Following this, Hughes, Keating and O'Connell~\cite{HKO1} used the characteristic polynomial to model the discrete moments of the derivative of the zeta function,
\begin{equation*}
J_k(T) := \frac{1}{N(T)} \sum_{0<\g_n\leq T} \left|\zeta'(\tfrac{1}{2}+\i\g_n)\right|^{2k} .
\end{equation*}
They calculated that for fixed $k$ subject to $\Re(k)>-3/2$,
\begin{equation*}
\E\left\{\frac{1}{N}\sum_{n=1}^N \left|Z_U'(\t_n)\right|^{2k} \right\} \sim \frac{G^2(k+2)}{G(2k+3)} N^{k(k+2)}
\end{equation*}
as $N\to\infty$, and they used this to conjecture the asymptotic form $J_k(T)$ should take for large $T$.
\begin{conj}\label{conj:HKO}
{\bf (Hughes, Keating and O'Connell).}
If all the zeros of the zeta function are all simple, then for fixed $k>-3/2$, as $T\to\infty$,
\begin{equation*}
\frac{1}{N(T)} \sum_{0<\g_n\leq T} \left|\zeta'(\tfrac{1}{2}+\i\g_n)\right|^{2k} \sim \frac{G^2(k+2)}{G(2k+3)} a(k) \left(\log \frac{T}{2\pi}\right)^{k(k+2)}
\end{equation*}
where $a(k)$ is given by \eqref{eq:a(k)}, and $G(\cdot)$ is the Barnes $G$--function.
\end{conj}
Again this conjecture is found to agree with all previously known results; when $k=-1$ (a conjecture of Gonek~\cite{Gon1}), when $k=0$ (trivial), and when $k=1$ (a theorem of Gonek~\cite{Gon0} under RH)\@. Also, extending a theorem due to Conrey, Ghosh and Gonek~\cite{CGG} (recorded in this paper as theorem~\ref{thm:CGG}) beyond its range of (proven) applicability, conjecture~\ref{conj:disc_deriv_k=2} states that $J_2(T) \sim  \frac{1}{1440\pi^2} \left(\log\frac{T}{2\pi}\right)^{8}$, which agrees perfectly with conjecture~\ref{conj:HKO}.

It is striking that conjecture~\ref{conj:KS} and conjecture
\ref{conj:HKO} have very similar form. The purpose of this
paper is to unify them as special cases of one result.

In the next section we will show
\begin{equation*}
\E_N \left\{ \left| Z_U\left(\t_1 + \frac{y}{N}\right)\right|^{2k}\right\} \sim \frac{G^2(k+1)}{G(2k+1)} F_k(y) N^{k^2}
\end{equation*}
where $F_k(y)$ is a certain function, independent of $N$, given in theorem~\ref{thm:displaced_mmt}.

We will then use this to conjecture that for $L=\frac{1}{2\pi}\log\frac{T}{2\pi}$,
\begin{equation*}
\frac{1}{N(T)}\sum_{0 < \g_n \leq T} \left| \zeta\left(\tfrac{1}{2}+\i \left(\g_n + \a/L\right) \right)\right|^{2k} \sim \frac{G^2(k+1)}{G(2k+1)} a(k) F_k(2\pi\a) \left(\log\frac{T}{2\pi}\right)^{k^2}
\end{equation*}
and will show, in section~\ref{sect:jnts_mmts_zeta}, that this conjecture contains conjectures~\ref{conj:KS} and~\ref{conj:HKO} as special cases ($\a\to\infty$ and $\a\to 0$ respectively). The conjecture is found to agree with a known result of Gonek~\cite{Gon0}, and the extension of the theorem due to Conrey, Ghosh and Gonek~\cite{CGG} cited above. These comparisons are discussed in section~\ref{sect:compare zeta function}.

\section{The random matrix calculation}

\begin{theorem}\label{thm:displaced_mmt}
For fixed $k$ with $\Re(k)>-1/2$, and for $x\leq AN$ with $A<\pi$ constant,
\begin{equation*}
\E_N \left\{ \left| Z_U\left(\t_1 + \frac{2x}{N}\right)\right|^{2k}\right\} = \frac{G^2(k+1)}{G(2k+1)} F_k(2x) N^{k^2} + \O\left(N^{k^2-1}\right)
\end{equation*}
where
\begin{equation}\label{eq:defn F_k(2x)}
F_k(2x) = x^2 j_k(x)^2 + x^2 j_{k-1}(x)^2 - 2k x j_k(x) j_{k-1}(x)
\end{equation}
where $j_n(x)$ are the spherical Bessel functions of the first kind.
\end{theorem}

\noindent\textbf{Proof.}
First note that
\begin{equation*}
\left|Z_U(\t_N+\b)\right|^{2k} = \left|1-e^{-\i\b}\right|^{2k} \prod_{n=1}^{N-1}\left|1-e^{\i(\t_n-\t_N-\b)}\right|^{2k}
\end{equation*}
The average of this over all $N\times N$ unitary matrices with Haar measure can be
written (\cite{Weyl,Mehta}) as an $N$--fold integral
\begin{multline*}
\E_N\left\{\prod_{n=1}^{N-1}\left|1-e^{\i(\t_n-\t_N-\b)}\right|^{2k} \right\} =\\
\frac{1}{N! (2\pi)^N} \int_{-\pi}^\pi \dots \int_{-\pi}^\pi
\prod_{1\leq i < j\leq N} \left|e^{\i\t_i} - e^{\i\t_j} \right|^2
\prod_{n=1}^{N-1}\left|1-e^{\i(\t_n-\t_N-\b)}\right|^{2k} \d\t_1\dots\d\t_N
\end{multline*}
Putting all the $j=N$ terms from the first product into the second, we
see that
\begin{equation}\label{eq:Z(t1+y)_factor}
\E_N\left\{ \left|Z_U\left(\t_N+\b\right)\right|^{2k}\right\} = \frac{1}{N} \left|2\sin\left(\tfrac{1}{2}\b\right)\right|^{2k} \E_{N-1}\left\{ \left|Z_{\widetilde U}(0)\right|^{2} \left|Z_{\widetilde U}(\b)\right|^{2k}\right\} ,
\end{equation}
where $Z_{\widetilde U}$ is the characteristic polynomial of an
$(N-1)\times(N-1)$ unitary matrix.

By rotation invariance of Haar measure,
\begin{equation*}
\E_{N-1}\left\{ \left|Z_{\widetilde U}(0)\right|^{2} \left|Z_{\widetilde U}(\b)\right|^{2k}\right\} = \E_{N-1}\left\{ \left|Z_{\widetilde U}(0)\right|^{2k} \left|Z_{\widetilde U}(\b)\right|^{2}\right\}
\end{equation*}
This is calculated in Theorem \ref{thm: jnt mmt Z(0)Z(b)} below, where
it is shown that for $\Re(k)>-1/2$,
\begin{multline*}
\E_N\left\{|Z_U(0)|^{2k} |Z_U(y/N)|^2 \right\} \\
\sim \frac{G^2(k+1)}{G(2k+1)} \sum_{p=0}^\infty \frac{k (k-1+p)! (k+p)!}{p! (2k+p)! (2k+1+2p)!} (-1)^p y^{2p} N^{(k+1)^2}
\end{multline*}
and substituting this into \eqref{eq:Z(t1+y)_factor} (where we put $\beta=y/N$) we see that
\begin{equation*}
\E_N\left\{ \left|Z_U\left(\t_N+\frac{y}{N}\right)\right|^{2k}\right\}
= \frac{G^2(k+1)}{G(2k+1)} F_k(y) N^{k^2}\left(1  +\O\left(\frac{1}{N}\right)\right)
\end{equation*}
where
\begin{equation} \label{eq:F_series}
F_k(y) =k \sum_{p=0}^\infty \frac{(k-1+p)! (k+p)!}{p! (2k+p)! (2k+1+2p)!}
(-1)^p y^{2k+2p}
\end{equation}

The spherical Bessel functions of the first kind are defined as
\begin{align*}
j_n(z) &= \sqrt{\frac{\pi}{2z}} J_{n+1/2}(z)\\
&= \sum_{m=0}^\infty \frac{(-1)^m (n+m)! }{m! (2n+2m+1)!2^{2m}} (2z)^{n+2m}
\end{align*}
where $J_\nu(z)$ is the $\nu$-th order Bessel function of the first kind. Hence,
\begin{equation*}
F_k(2x) = x^2 j_k(x)^2 + x^2 j_{k-1}(x)^2 - 2k x j_k(x) j_{k-1}(x)
\end{equation*}
which can be seen by comparing the Taylor expansions.

The above large-$N$ asymptotics are for $x=o(N)$. This can be extended to $|x| \leq A N$ for $A<\pi$ an arbitrary constant as follows: Let $\beta$ be a fixed constant subject to $0<\beta<2\pi$. By \eqref{eq:Z(t1+y)_factor} and the results of Basor~\cite{B},
\begin{align*}
\E_N \left\{ \left| Z_U\left(\t_1 + \beta\right)\right|^{2k}\right\} &= \frac{1}{N} \left|2\sin\left(\tfrac{1}{2}\beta\right)\right|^{2k} \E_{N-1}\left\{ \left|Z_{\widetilde U}(0)\right|^{2k} \left|Z_{\widetilde U}(\beta)\right|^{2}\right\}\\
&\sim \frac{G^2(k+1)}{G(2k+1)} N^{k^2} .
\end{align*}
If one lets $x=N\beta/2$ then the large-$x$ asymptotics of Bessel functions (see, for example, chapter 9 of \cite{AbramS}) implies that
\begin{equation*}
F_k(N\beta) = 1 +\O\left(\frac{1}{N}\right) .
\end{equation*}
and so theorem~\ref{thm:displaced_mmt} gives the correct first order term as $N\to\infty$, for $|x| \leq A N$ with $A<\pi$ a constant.
\qed

\noindent\textbf{Remark.}
When $n$ is an integer,
\begin{equation*}
j_n(x) = (-1)^n x^n \left( \frac{1}{x}\frac{\d}{\d x}\right)^n \frac{\sin x}{x}
\end{equation*}
which leads to a neat evaluation of $F_k(2x)$ for integer $k$, the first few being:
\begin{gather}
F_1(2x) = \frac{x^2 - \sin^2(x)}{x^2}\nonumber\\
F_2(2x) = \frac{x^4-3x^2 + 3x\sin(2x) + (2x^2-3)\sin^2(x)}{x^4}\label{eq:F_2(2x)}\\
F_3(2x) = \frac{x^6-3x^4-45x^2 + (-12x^3+45x)\sin(2x) + (-3x^4+72x^2-45)\sin^2(x)}{x^6}\nonumber
\end{gather}

\begin{theorem}\label{thm: jnt mmt Z(0)Z(b)}
\begin{multline*}
\E_N\left\{ |Z(0)|^{2k} |Z(\b)|^2 \right\} =  M_N(2k) N! (N+2k)! \times\\
\times \sum_{n=0}^{N}  \frac{(2\sin(\b/2))^{2n}}{n!(2k+n)!} \sum_{m=0}^{N-n} \frac{(N+k-m)!(k+m+n)!}{(N-n-m)!m!}  e^{\i\b(2m-N+n)}
\end{multline*}
where $M_n(2k)$ is given in \eqref{eq:defn M_N(2k)}. If $\frac{y}{N}\to 0$ as $N\to\infty$, then for $\Re(k)>-1/2$,
\begin{multline*}
\E_N\left\{ |Z(0)|^{2k} |Z(\frac{y}{N})|^2 \right\} = \frac{G^2(k+1)}{G(2k+1)}\sum_{p=0}^\infty \frac{k (k-1+p)! (k+p)!}{p! (2k+p)! (2k+1+2p)!} (-1)^p y^{2p} \\
\times N^{(k+1)^2} \left(1+\O(\frac{1}{N})\right)
\end{multline*}
\end{theorem}

\noindent \textbf{Proof.}
If $k$ in an integer, then
\begin{multline*}
\E_N\left\{\prod_{n=1}^N |1-e^{\i\t_n}|^{2k} |1-e^{\i\t_n}e^{-\i\b}|^2\right\}\\
 = \frac{1}{(2\pi)^N N!} \idotsint_{-\pi}^\pi \prod_{1\leq j<m\leq N} (e^{\i\t_j}-e^{\i\t_m})(e^{-\i\t_j}-e^{-\i\t_m}) \\
\times \prod_{n=1}^N (1-e^{\i\t_n})^k(1-e^{-\i\t_n})^k (1-e^{\i\t_n}e^{-\i\b})(1-e^{-\i\t_n}e^{\i\b}) \;\d\t_n
\end{multline*}
which equals, after some simple manipulation of the terms
\begin{multline*}
\frac{(-1)^{N(N-1)/2+kN+N}e^{-\i N\b}}{(2\pi)^N N!}\times\\
\idotsint_{-\pi}^\pi \prod_{1\leq j<m\leq N}
(e^{\i\t_j}-e^{\i\t_m})^2 \prod_{n=1}^N (e^{\i\t_n})^{-N-k}
(1-e^{\i\t_n})^{2k} (e^{\i\t_n}-e^{\i\b})^2 \;\d\t_n
\end{multline*}
Therefore,
\begin{multline*}
\E_N\left\{ |Z_U(0)|^{2k} |Z_U(\b)|^2 \right\} = \frac{e^{-\i N\b}}{N!} (-1)^{N(N-1)/2+kN+N}\\
\times \CT \left\{\prod_{1\leq j<m\leq N} (t_j-t_m)^2
\prod_{n=1}^N \frac{1}{t_n^{N+k}} (1-t_n)^{2k}
(t_n-e^{\i\b})^2\right\}
\end{multline*}
where $\CT\{\cdot\}$ denotes the constant term in the Laurent expansion in the variables $t_1,\dots,t_N$. The constant term equals (by lemma 1 of \cite{Forr})
\begin{equation*}
\lim_{y\to 0} y^N \idotsint_0^1 \prod_{1\leq j<m\leq N}
(t_j-t_m)^2 \prod_{n=1}^N t_n^{-N-k+y-1} (1-t_n)^{2k}
(t_n-e^{\i\b})^2 \;\d t_n
\end{equation*}

Kaneko \cite{Kaneko} has evaluated this integral (which is a generalization of Selberg's integral) as
\begin{multline*}
\prod_{j=1}^N \frac{\Gamma(1+j) \Gamma(j+y-N-k)\Gamma(j+2k+1)}{\Gamma(j+y+k+1)}\\
\sum_{m,n=0}^\infty \frac{(-N)_{m+n} (y+k+1)_{m+n}}{(y-N-k)_m (2k+1)_n} \frac{e^{2m\i\b} (1-e^{\i\b})^{2n}}{m! n!}
\end{multline*}
where $(a)_n = a(a+1)\dots(a+n-1)=\Gamma(a+n)/\Gamma(a)$.
Since we have assumed that $k$ is an integer,
\begin{equation*}
\lim_{y\to 0} y \Gamma(y+j-N-k) =\frac{(-1)^{N+k-j}}{\Gamma(N+k-j+1)}
\end{equation*}
and so we have
\begin{multline*}
\E_N\left\{ |Z_U(0)|^{2k} |Z_U(\b)|^2 \right\} = \frac{e^{-\i N\b}}{N!} \prod_{j=1}^N \frac{\Gamma(1+j) \Gamma(j+2k+1)}{\Gamma(j+k+1)\Gamma(N+k-j+1)}\\
\times \sum_{m,n=0}^\infty \frac{(-N)_{m+n} (k+1)_{m+n}}{(-N-k)_m (2k+1)_n} \frac{e^{2m\i\b} (1-e^{\i\b})^{2n}}{m! n!}
\end{multline*}
Expanding everything out in terms of the gamma function:
\begin{multline*}
\E_N\left\{ |Z_U(0)|^{2k} |Z_U(\b)|^2 \right\} =  M_N(2k)
\frac{\Gamma(N+1+2k) N!}{\Gamma(N+1+k)^2} \\
\times \sum_{n=0}^{N}  \frac{(2\sin(\b/2))^{2n}}{n!\Gamma(2k+n+1)} \sum_{m=0}^{N-n} \frac{\Gamma(N+k+1-m)\Gamma(k+n+1+m)}{\Gamma(N+1-n-m)!}  \frac{e^{\i\b(2m-N+n)} }{m!}
\end{multline*}
where $M_N(2k)$ is defined in \eqref{eq:defn M_N(2k)}. Observe
that the inner summand is invariant as $m \longrightarrow N-n-m$,
and so the inner sum is in fact a sum of cosines (and thus the
series expansion in $\b$ contains only even powers of $\b$).
Furthermore observe that both sides of the equation are analytic
functions of $k$ (for $\Re(k)\geq 0$), both sides can be easily
bounded bounded by $\O(2^{2N\Re(k)})$ (for large $k$, with $N$ and
$\beta$ fixed), and the two sides are equal at the positive
integers. Thus by Carlson's theorem (see \S17 of \cite{Mehta}, for
example), the restriction on $k$ being an integer is no longer
required, and the left-hand side equals the right-hand side for
all complex $k$.

Now, using the fact that $(-a)_n = (-1)^n \Gamma(a+1)/\Gamma(a+1-n)$, we have
\begin{multline*}
\sum_{m=0}^{N-n} \frac{\Gamma(N+k+1-m)\Gamma(k+n+1+m)}{\Gamma(N+1-n-m)}  \frac{e^{\i\b(2m-N+n)} }{m!}\\
=\frac{\Gamma(N+k+1) \Gamma(k+n+1)}{\Gamma(N-n+1)} \Hyper{-N+n,k+n+1}{-N-k}{e^{2\i\b}}e^{\i(n-N)\b}
\end{multline*}
where
\begin{align*}
\Hyper{a,b}{c}{z} &= \sum_{m=0}^\infty \frac{(a)_m (b)_m}{(c)_m}\frac{z^m}{m!}
\end{align*}
is the Gauss hypergeometric function. (If $a,c$ are negative integers with $-a<-c$ then it is a polynomial of degree $-a$).

Applying the quadratic hypergeometric transformation 15.3.26 of \cite{AbramS} we get
\begin{multline}\label{eq:hyper_quad_trans}
\Hyper{-N+n,k+n+1}{-N-k}{e^{2\i\b}} e^{\i(n-N)\b} \\
=  \Hyper{-\tfrac{1}{2}N+\tfrac{1}{2}n,-\tfrac{1}{2}N+\tfrac{1}{2}n+\tfrac{1}{2} }{ -N-k }{\frac{1}{\cos^2(\b)}}(2\cos\b)^{N-n}
\end{multline}
For $m$ a positive integer,
\begin{equation*}
\Hyper{-m,b}{c}{z} = \frac{\Gamma(1-c-m)\Gamma(1-c+b)}{\Gamma(1-c)\Gamma(1-c-m+b)} \Hyper{-m,b}{1-c-m+b}{1-z}
\end{equation*}
and so we see that the right-hand side of \eqref{eq:hyper_quad_trans} equals
\begin{multline*}
\frac{\Gamma(\tfrac{1}{2}N+\tfrac{1}{2}n+k+1) \Gamma(\tfrac{1}{2}N+\tfrac{1}{2}n+k+\tfrac{3}{2})}{\Gamma(N+k+1) \Gamma(k+n+\tfrac{3}{2})}(2\cos\b)^{N-n} \times\\
\times \Hyper{-\tfrac{1}{2}N+\tfrac{1}{2}n,-\tfrac{1}{2}N+\tfrac{1}{2}n+\tfrac{1}{2} }{ n+k+\tfrac{3}{2} }{ 1-\frac{1}{\cos^2(\b)}}
\end{multline*}

Therefore, we have proven
\begin{multline*}
\E_N\left\{ |Z_U(0)|^{2k} |Z_U(\b)|^2 \right\} = \\
\sum_{n=0}^{N}  \sum_{m=0}^{\lfloor\frac12(N-n)\rfloor} T(N,k,m,n) (-1)^m  (2\sin(\b/2))^{2n} (\sin\b)^{2m} (\cos\b)^{N-n-2m}
\end{multline*}
where
\begin{multline*}
T(N,k,m,n) =  M_N(2k)  \frac{N!(N+2k)! (\tfrac{1}{2}N+\tfrac{1}{2}n+k)! (\tfrac{1}{2}N+\tfrac{1}{2}n+k+\tfrac{1}{2})!}{(N-n)! (N+k)!^2} \\
\times  \frac{2^{N-n} (-\tfrac12N+\tfrac12n)_m (-\tfrac12N+\tfrac12n+\tfrac12)_m (k+n)!}{m! n!(2k+n)!(n+k+\tfrac12+m)!}
\end{multline*}
Observe that for fixed $k,m,n$ with $\Re(k)>-1/2$,
\begin{multline*}
T(N,k,m,n) = \frac{G^2(k+1)}{G(2k+1)} \frac{(k+n)! (k+n+m)! }{m! n!(2k+n)! (2k+2n+2m+1)!} \\
\times N^{(k+1)^2 +2n+2m}\left(1+\O(\frac{1}{N})\right)
\end{multline*}
and so
\begin{multline*}
\E_N\left\{ |Z(0)|^{2k} |Z(y/N)|^2 \right\} \sim \frac{G^2(k+1)}{G(2k+1)} N^{(k+1)^2}\\
\times \sum_{n=0}^{\infty} \sum_{m=0}^\infty  (-1)^m \frac{(k+n)! (k+n+m)! }{m! n!(2k+n)! (2k+2n+2m+1)!} y^{2m+2n}
\end{multline*}
in the sense that for each fixed integer $h$, the coefficient of $y^h$ on the Taylor expansion of the left-hand side converges to that of the right-hand side as $N\to\infty$.

Finally we show that
\begin{multline*}
\sum_{n=0}^{\infty} \sum_{m=0}^\infty  (-1)^m \frac{(k+n)! (k+n+m)! }{m! n!(2k+n)! (2k+2n+2m+1)!} y^{2m+2n}\\
= \sum_{p=0}^\infty \frac{k (k-1+p)! (k+p)!}{p! (2k+p)! (2k+1+2p)!}
(-1)^p y^{2p}
\end{multline*}
This can be proven by comparing the coefficients of $y^{2p}$. That is, we wish to show for all integer $p\geq 0$,
\begin{equation*}
\sum_{n=0}^p \frac{(-1)^{p-n} (k+n)! (k+p)! }{(p-n)! n!(2k+n)! (2k+2p+1)!}
= \frac{(-1)^p k (k-1+p)! (k+p)!}{p! (2k+p)! (2k+1+2p)!}
\end{equation*}
This is equivalent to showing
\begin{equation}\label{eq:WZ_sum}
\sum_{n=0}^p \frac{(-1)^n (k+n)! p! (2k+p)!}{k (p-n)! n!(2k+n)! (k-1+p)!} = 1
\end{equation}
for all integer $p\geq 0$, and this we shall do by creating the Wilf-Zeilberger pair \cite{A=B}: Denote the summand in \eqref{eq:WZ_sum} by $F(p,n)$, and observe that
\begin{equation}\label{eq:WZ_sum this}
F(p+1,n)-F(p,n) = G(p,n+1)-G(p,n)
\end{equation}
where
\begin{equation*}
G(p,n) = \frac{(2k+n) n}{(n-p-1) (k+p)} F(p,n)
\end{equation*}
(the $(2k+n) n/(n-p-1) (k+p)$ being calculated by Zeilberger's
algorithm). Summing both sides of \eqref{eq:WZ_sum this}
over $n$, we see that the right-hand side telescopes to zero, which
shows that the left-hand side of \eqref{eq:WZ_sum} must be a constant,
independent of $p$. Putting $p=0$, direct calculation shows that constant is $1$.
\qed

\section{Conjecture about the zeta function}
\label{sect:jnts_mmts_zeta}

\begin{conj}\label{conj:jnt_mmts}
For fixed $k$ subject to $\Re(k)>-1/2$,
\begin{equation*}
\frac{1}{N(T)} \sum_{0 < \g_n \leq T} \left| \zeta\left(\tfrac{1}{2}+\i \left(\g_n + \a/L\right) \right)\right|^{2k} \sim \frac{G^2(k+1)}{G(2k+1)} a(k) F_k(2\pi\a) \left(\log\frac{T}{2\pi}\right)^{k^2}
\end{equation*}
as $T\to\infty$, uniformly in $\a$ for $|\a| \leq L$, where $L=\frac{1}{2\pi}\log\frac{T}{2\pi}$ is the density of
zeros of height $T$. $G(\cdot)$ is the Barnes $G$--function, $a(k)$ is
given by \eqref{eq:a(k)}, and $F_k(2\pi\a)$ is given in theorem~\ref{thm:displaced_mmt}.
\end{conj}

If this conjecture is true, then we are able to prove conjecture~\ref{conj:HKO} and a variant of the Keating-Snaith conjecture (conjecture~\ref{conj:KS}).
\begin{cor}\label{conj:disc_deriv_from_jnt_mmts}
If conjecture~\ref{conj:jnt_mmts} is true, then
\begin{equation*}
\frac{1}{N(T)} \sum_{0 < \g_n \leq T} \left| \zeta'\left(\tfrac{1}{2}+\i \g_n \right)\right|^{2k} \sim \frac{G^2(k+2)}{G(2k+3)} a(k) \left(\log\frac{T}{2\pi}\right)^{k(k+2)}
\end{equation*}
\end{cor}

\noindent\textbf{Proof.}
By the definition of differentiation,
\begin{equation*}
\left| \zeta'\left(\tfrac{1}{2}+\i \g_n \right)\right|^{2k} = L^{2k} \lim_{a\to 0} \frac{\left| \zeta\left(\tfrac{1}{2}+\i \left(\g_n + \frac{\a}{L}\right)\right)\right|^{2k}}{\a^{2k}}
\end{equation*}
From \eqref{eq:F_series} we have
\begin{equation*}
\lim_{\a\to 0} \frac{F(2\pi\a)}{\a^{2k}} = (2\pi)^{2k} \frac{k! k!}{(2k)! (2k+1)!}
\end{equation*}
so applying conjecture~\ref{conj:jnt_mmts} and using uniformity to swap the $\a\to 0$ and $N\to\infty$ limits, we have
\begin{align*}
\frac{1}{N(T)} \sum_{0 < \g_n \leq T} \left| \zeta'\left(\tfrac{1}{2}+\i \g_n \right)\right|^{2k} &\sim \frac{G^2(k+1)}{G(2k+1)} a(k) L^{2k}  \frac{(2\pi)^{2k} k! k!}{(2k)! (2k+1)!} \left(\log\frac{T}{2\pi}\right)^{k^2}\\
&= \frac{G^2(k+2)}{G(2k+3)} a(k) \left(\log\frac{T}{2\pi}\right)^{k(k+2)}
\end{align*}
as required.
\hfill$\Box$

\begin{cor}\label{conj:variant KS}
From conjecture~\ref{conj:jnt_mmts} it follows that for $\beta>0$ fixed
\begin{equation}\label{eq:variant_KS}
\frac{1}{N(T)} \sum_{0 < \g_n \leq T} \left| \zeta\left(\tfrac{1}{2}+\i (\g_n+\beta) \right)\right|^{2k} \sim \frac{G^2(k+1)}{G(2k+1)} a(k) \left(\log\frac{T}{2\pi}\right)^{k^2}
\end{equation}
\end{cor}
\noindent\textbf{Proof.}
The asymptotics as $z\to\infty$ of the spherical Bessel function (see chapters 9 and 10 of~\cite{AbramS}) are
\begin{equation*}
j_n(z) \sim
\begin{cases}
\frac{1}{z}(-1)^{(n+1)/2} \cos z & \text{ if $n$ is odd}\\
\frac{1}{z}(-1)^{n/2} \sin z & \text{ if $n$ is even}
\end{cases}
\end{equation*}
and putting this into \eqref{eq:defn F_k(2x)} we have
\begin{equation*}
\lim_{\a\to\infty} F_k(2\pi\a) = 1
\end{equation*}
Setting $\a=L\beta$ in conjecture~\ref{conj:jnt_mmts} completes the proof.
\qed

\noindent\textbf{Remark.} Note that corollary~\ref{conj:variant
KS} can be thought of as a variant of conjecture \ref{conj:KS}.
This is because one expects the mean of $|\zeta(1/2+\i t)|^{2k}$
to be independent of the average taken, and the Keating-Snaith
conjecture is a result about the continuous mean, whereas
corollary~\ref{conj:variant KS} is a result about a discrete mean.
To see this, recall that the zeros get denser higher up the
critical line, and so if $\beta>0$ is fixed and $\gamma_n$ is
random, one might expect $\zeta\left(\tfrac{1}{2}+\i (\g_n+\beta)
\right)$ to be random (whereas, if $\beta$ was small, it would be
highly influenced by the fact that $\zeta\left(1/2+\i \g_n \right)
= 0$). The left-hand side of (\ref{eq:variant_KS}) averages this,
and thus acts like a discrete mean of $|\zeta(1/2+\i t)|^{2k}$.

\subsection{Comparison with the zeta function}\label{sect:compare zeta function}

Gonek~\cite{Gon0} showed that if the Riemann Hypothesis is true then
\begin{equation*}
\frac{1}{N(T)} \sum_{0 < \g_n \leq T} \left| \zeta\left(\tfrac{1}{2}+\i \left(\g_n + \a/L\right) \right)\right|^{2} \sim \left(1-\left(\frac{\sin(\pi\a)}{\pi\a}\right)^2\right) \log \frac{T}{2\pi}
\end{equation*}
uniformly in $\a$ for $|\a| \leq L/2$, which is in perfect agreement with conjecture~\ref{conj:jnt_mmts} when $k=1$.

There is no proof of the conjecture for $k=2$ (unlike conjecture~\ref{conj:KS} which is proven for $k=1$ and $2$).
But there are theorems along the lines of conjecture~\ref{conj:jnt_mmts} for $k=2$:
\begin{theorem}\label{thm:CGG}
{\bf (Conrey, Ghosh and Gonek~\cite{CGG}.)}
Assume GRH and let $A(s) = \sum_{n\leq x} n^{-s}$ where $x = \left(\frac{T}{2\pi}\right)^\eta$ for some $\eta\in(0,\tfrac{1}{2})$. Then,
\begin{multline*}
\frac{1}{N(T)}\sum_{0 < \g_n \leq T} \left|\zeta A \left(\tfrac{1}{2}+\i\left(\g_n+\a/L\right)\right)\right|^{2} \sim \frac{6}{\pi^2} \sum_{j=0}^\infty \frac{(-1)^{j+1} (2\pi\a)^{2j+2}}{(2j+5)!} \times\\
\times \left(-\eta^2+\tfrac{1}{3}(2j+5)\eta^3-\frac{2j+5}{j+3}\eta^{2j+6}+\eta^{2j+7}+\eta^2(1-\eta)^{2j+5}\right)\left(\log\frac{T}{2\pi}\right)^4
\end{multline*}
uniformly for bounded $\a$.
\end{theorem}
\noindent(We have slightly changed notation from~\cite{CGG}, to be consistent with our definition of $L=\frac{1}{2\pi}\log\frac{T}{2\pi}$).

Putting $\eta=1$ in the above (which, as it stands, is not allowed under the conditions of the theorem) then $A(\tfrac{1}{2}+\i t) = \zeta(\tfrac{1}{2}+\i t)+\O(t^{-1/2})$, and we have
\begin{multline*}
\frac{1}{N(T)}\sum_{0 < \g_n \leq T} \left|\zeta^2\left(\tfrac{1}{2}+\i\left(\g_n+\a/L\right)\right)\right|^{2} \\
\sim \frac{4}{\pi^2} \sum_{j=1}^\infty \frac{(-1)^{j+1} (2\pi\a)^{2j+2} }{(2j+6)!} (2j^2+5j) \left(\log \frac{T}{2\pi}\right)^4
\end{multline*}
Note that
\begin{multline*}
\frac{4}{\pi^2} \sum_{j=1}^\infty \frac{(-1)^{j+1} (2\pi\a)^{2j+2}}{(2j+6)!} (2j^2+5j) =\\
\frac{1}{12}a(2)\frac{(2\pi^2\a^2-3)\sin^2(\pi\a) + 3\pi\a\sin(2\pi\a)+(\pi\a)^4-3(\pi\a)^2}{(\pi\a)^4}
\end{multline*}
which is what is predicted in conjecture~\ref{conj:jnt_mmts}.

That is, from a purely number theoretical calculation involving no random matrix theory, we have
\begin{conj}\label{conj:jnt_mmts_k=2}
Assuming that $\eta=1$ is permissible in theorem~\ref{thm:CGG} then
\begin{equation*}
\frac{1}{N(T)}\sum_{0 < \g_n \leq T} \left|\zeta\left(\tfrac{1}{2}+\i\left(\g_n+\a/L\right)\right)\right|^{4} \sim \frac{1}{2\pi^2} F_2(2\pi\a) \left(\log \frac{T}{2\pi}\right)^4
\end{equation*}
where $F_2(2x)$ is given in \eqref{eq:F_2(2x)}.
\end{conj}

So, following the proof of corollary~\ref{conj:disc_deriv_from_jnt_mmts}, we may deduce
\begin{cor}\label{conj:disc_deriv_k=2}
If conjecture~\ref{conj:jnt_mmts_k=2} is true then
\begin{equation*}
\frac{1}{N(T)} \sum_{0 < \g_n \leq T} \left| \zeta'\left(\tfrac{1}{2}+\i \g_n \right)\right|^{4} \sim \frac{1}{1440\pi^2} \left(\log\frac{T}{2\pi}\right)^{8} .
\end{equation*}
\end{cor}
Note that this is the same answer that one gets from putting $k=2$ into conjecture~\ref{conj:HKO}.

\section{Acknowledgments}

Early versions of the results described in this paper formed part of my PhD thesis~\cite{hughes}, and I would like to thank my advisors, Jon Keating and Neil O'Connell, for their help and advice. Part of this research was carried out while I was supported by the EC TMR network ``Mathematical aspects of Quantum Chaos'', EC-contract no HPRN-CT-2000-00103.

\end{document}